\documentclass[12pt]{amsart}
\usepackage{latexsym, amssymb, amsmath, amscd}
 \usepackage{eucal}

\theoremstyle{plain}
    \newtheorem{rema}{Remark}[section]
    \newtheorem{propo}[rema]{Proposition}

   \newtheorem{theo}[rema]{Theorem}
 \newtheorem{conj}[rema]{Conjecture}
   
    \newtheorem{lemma}[rema]{Lemma}
    \newtheorem{corol}[rema]{Corollary}
     \newtheorem{exam}[rema]{Example}
  \newtheorem{rmk}[rema]{Remark}

	\newcommand{\nno}{\nonumber}

	\newcommand{\p}{\partial}

 \newcommand{\pf}{{\it Proof:}\hspace{2ex}}
 
 \newcommand{\epfv}{\hspace{1em}$\Box$\vspace{1em}}

\newcommand{\bZ}{{\mathbb Z}}

\newcommand{\bN}{{\mathbb N}}

\newcommand{\bC}{\mathbb C}

\newcommand{\cR}{{\mathcal R}}
\newcommand{\cL}{{\mathcal L}}

\newcommand{\cM}{{\mathcal M}}

\newcommand{\cJ}{{\mathcal J}}

\newcommand{\cA}{{\mathcal A}}

\newcommand{\cD}{{\mathcal D}}

\newcommand{\xiz}{{\xi, z}}

\newcommand{\lin}{1\le i\le n}

\renewcommand{\theequation}{\thesection.\arabic{equation}}
\renewcommand{\therema}{\thesection.\arabic{rema}}
\newcommand{\EAn}{\end{align*}}

\title[Abhyankar-Gurjar's Formula and the Jacobian Conjecture]
{New Proofs for the Abhyankar-Gurjar Inversion Formula and the Equivalence of the Jacobian Conjecture and the Vanishing Conjecture}

  \author{Wenhua Zhao}      
    \date{\today}

\begin{document}

\begin{abstract}
We first give a new proof and also a new formulation for the Abhyankar-Gurjar inversion formula for formal maps of affine spaces. We then use the reformulated 
Abhyankar-Gurjar formula to give a more straightforward proof for the equivalence of the {\it Jacobian conjecture} with a special case of the 
{\it vanishing conjecture} of (homogeneous) quadratic differential operators with constant coefficients.  
\end{abstract}

   

\thanks{The author has been partially supported 
by NSA Grant H98230-10-1-0168}

 \bibliographystyle{alpha}
    \maketitle


\renewcommand{\theequation}{\thesection.\arabic{equation}}
\renewcommand{\therema}{\thesection.\arabic{rema}}
\setcounter{equation}{0}
\setcounter{rema}{0}
\setcounter{section}{0}

\section{\bf Introduction} \label{S1}

Let $K$ be a field of characteristic zero and 
$z=(z_1, z_2, \dots, z_n)$ $n$ 
commutative free variables. We denote by $\cA_K[z]$ 
(resp., $\cA_K[[z]]$) the algebra over $K$ 
of polynomials (resp., formal power series) 
in $z$. Set $\p_i\!:=\p_{z_i}$ $(\lin)$ 
and $\p\!:=(\p_1, \p_2, ..., \p_n)$.

Recall that the {\it Jacobian conjecture}  
posed by O. H. Keller \cite{Ke} in 1939  
claims that {\it any polynomial map 
$F=(F_1, F_2, ..., F_n)$ 
of $K^n$ with Jacobian 
$\cJ F(z)\!:=\det ( \p_j F_i )\equiv 1$ 
must be an automorphism of $K^n$}. 
Despite intense study 
from mathematicians in the last seventy years, 
the conjecture is still open
even for the case $n=2$. In 1998, 
S. Smale \cite{Sm} included the  {\it Jacobian conjecture}  
in his list of $18$ fundamental mathematical problems 
for the $21$st century. For more history and 
known results on the  {\it Jacobian conjecture}, 
see \cite{BCW}, \cite{E} and 
references therein.

In the study of the {\it Jacobian conjecture}, the following 
two remarkable reductions play very important roles. 

The first one is the so-called {\it homogeneous reduction}, which 
was achieved independently by Bass, Connell and Wright \cite{BCW} and 
Jag\v zev \cite{J}. The reduction says that in order to prove 
or disprove the {\it Jacobian conjecture}, it suffices to study of polynomial 
maps of the form $F=z-H(z)$ with each component $H_i(z)$ of 
$H(z)\in \cA_K[z]^{\times n}$ being zero or 
homogeneous of degree three. 

The second reduction is the so-called {\it symmetric reduction},  
which was achieved independently by de Bondt and 
van den Essen \cite{BE} and 
Meng \cite{Me}. This reduction says that 
the Jacobian matrix of the polynomial maps 
of the form above may be further assumed 
to be symmetric.   

Based on the two reductions above and  
some results in \cite{BurgersEq} on a  
deformation of polynomial maps, the author 
showed in \cite{HNP} and \cite{GVC} that  
the  {\it Jacobian conjecture}  is actually 
equivalent to the following {\it vanishing conjecture} 
on quadratic differential operators of $\cA_K[z]$, 
i.e.,  the differential operators which can be 
written (uniquely) as quadratic forms in 
$\p$ over $K$.

\begin{conj} \label{VC} $(${\bf The Vanishing Conjecture}$)$
Let $P(z)\in \cA_K[z]$ and $\Lambda$ a quadratic 
differential operator of $\cA_K[z]$ with 
$\Lambda^m(P^m)=0$ for all $m\ge 1$.  
Then we have $\Lambda^m (P^{m+1})=0$ when $m\gg 0$.
\end{conj}

Actually, as shown in \cite{GVC}, 
the  {\it Jacobian conjecture} will follow 
if one can just show that Conjecture \ref{VC} 
holds for a (single) sequence of quadratic differential 
operators $\{\Lambda_n\,|\, n\in \bN\}$ such that the rank  
${\rm rk}\,\Lambda_n\to \infty$ 
as $n\to \infty$,  
where the {\it rank} 
${\rm rk}\,\Lambda_n$ 
of quadratic differential operators  
is defined to be the rank  
of the corresponding quadratic 
forms. For the study of some more general forms 
of Conjecture \ref{VC}, see \cite{GVC}, 
\cite{EZ}, \cite{EWZ} and \cite{IC}. 

Another important result on the {\it Jacobian conjecture}  
is the Abhyankar-Gurjar inversion formula, which provides 
a nice formula for the formal inverse maps 
of polynomial maps (see Theorem \ref{AGF}).  
This formula was first proved by Gurjar (unpublished) 
and later by Abhyankar \cite{Ab} in a simplified form   
(see also \cite{BCW}). 

One remark on the Abhyankar-Gurjar inversion formula is as follows. 
By comparing Eq.\,(\ref{AGF-e1}) with Eq.\,(\ref{AGF-pe2}), 
we see that up to the factor $\cJ F$ (the Jacobian of 
the polynomial map $F$), the differential operator in the  
Abhyankar-Gurjar inversion formula coincides with 
the one obtained by applying the well-known 
anti-involution $\tau$ (see Eq.\,(\ref{Def-tau})) 
of the Weyl algebra to the differential operator in the 
Taylor series expansion of $u(F)=u(z-H)$ at $z$ 
(see Eqs.\,(\ref{AGF-pe2})). 
In particular, for polynomial 
maps $F$ that satisfies the Jacobian condition 
$\cJ F\equiv 1$, these two differential operators 
becomes identical! But, unfortunately, 
from the proof given in \cite{Ab} (see also \cite{BCW}), 
which to the author's best knowledge is also 
the only published proof in the literature,    
it is not easy to see why such a resemblance 
or connection exists at the first place.   

In this paper,  we first use some results 
in \cite{Exp} to give a new proof for 
the Abhyankar-Gurjar inversion formula, which 
provides a good explanation for the resemblance or 
connection mentioned above. We then use a 
result in \cite{T-Deform} to 
give a new formulation for the Abhyankar-Gurjar 
inversion formula, and also derive some of its
consequences. Finally, by using this re-formulated  
Abhyankar-Gurjar inversion formula  
and some of its consequences derived 
in this paper, we give a new proof for 
the equivalence of  the  {\it Jacobian conjecture}  
with a special case of {\it vanishing conjecture}, 
Conjecture \ref{VC} (see  Theorem \ref{JC=VC}).

Comparing with the proof of the equivalence of the 
 {\it Jacobian conjecture}  
and Conjecture \ref{VC} given in 
\cite{HNP} and \cite{GVC}, the proof given  
here is more straightforward. But, contrast to 
the proof in \cite{HNP} and \cite{GVC},  
it does not show the equivalence of 
the  {\it Jacobian conjecture}  with Conjecture 
\ref{VC} itself but only with 
a special case of it.

The arrangement of this paper is as follows.

In Section \ref{S2}, we discuss some $\bZ$-gradings and 
the associated completions  
of polynomial algebras and the Weyl algebras.  
The main aims of this section  
are, first, to lay down a firm ground setting for  
the arguments in the rest of paper and, second,  
to show that certain results on polynomial algebras and 
the Weyl algebras can be extended to 
their completions constructed in this section. 
These extended results such as 
Proposition  \ref{Phi-Trans} and 
Lemma \ref{tau-lemma} will 
play crucial roles in our 
later arguments. 

In Section \ref{S3}, we first give a new proof 
and also a new formulation for the Abhyankar-Gurjar 
inversion formula. We then derive some consequences 
from the re-formulated Abhyankar-Gurjar 
inversion formula, which will be needed  
in the next section.

In Section \ref{S4}, we give a new proof for the equivalence 
of the  {\it Jacobian conjecture}  with a special case of the 
 {\it vanishing conjecture}. In addition, a special case 
of a conjecture proposed in \cite{HNP} is also 
proved here (see  Conjecture \ref{Stabe-Nil} 
and Proposition \ref{Lambda-Nil-2}).

Finally, one remark on the base field of this paper 
is as follows.

Unless stated otherwise, the base field 
$K$ throughout the paper is assumed 
to be a field of characteristic zero (instead of  
the complex number field $\bC$). 
This more general choice of   
the base field is necessary, 
or at least convenient, 
for some base field extension  
arguments in the paper.  
Note also that some of the 
results quoted in this paper were only 
proved over the complex number field 
$\bC$ in the quoted references. 
But, either by applying Lefschetz's principle 
or by going through their proofs 
(which actually work equally well for 
all the fields of characteristic zero),  
the reader can easily see 
that these quoted results  
also hold as long as the base field 
has characteristic zero.  
The same remark will not 
be repeated for each of 
these quoted results in the 
context of this paper.

\renewcommand{\theequation}{\thesection.\arabic{equation}}
\renewcommand{\therema}{\thesection.\arabic{rema}}
\setcounter{equation}{0}
\setcounter{rema}{0}

\section{\bf Some $\bZ$-gradings and Associated Completions of Polynomial Algebras 
and the Weyl Algebras} \label{S2}

Let $K$ be any field of characteristic zero and 
$\xi=(\xi_1,  ..., \xi_n)$ and 
$z=(z_1,  ..., z_n)$ $2n$ 
free commutative variables. 
We denote by $\cA_K[\xi]$, 
$\cA_K[z]$ and $\cA_K[\xiz]$  
the polynomial algebras over $K$
in $\xi$, $z$ and $(\xi, z)$, 
respectively. 
The corresponding formal power series algebras 
will be denoted respectively by 
$\cA_K[[\xi]]$, $\cA_K[[z]]$ 
and $\cA_K[[\xiz]]$. 
We denote by $\cD_K[z]$  
the differential operator algebra 
or the Weyl algebra of 
the polynomial algebra 
$\cA_K[z]$.

In this section, we consider some 
$\bZ$-gradings on the polynomial 
algebra $\cA_K[\xiz]$ and 
the Weyl algebra  $\cD_K[z]$, 
and also the associated 
completions of $\cA_K[\xiz]$ 
and $\cD_K[z]$ with respected 
to these $\bZ$-gradings.
All formal power series and 
differential operators involved 
in this paper will lie in the completions constructed 
in this section for $\cA_K[\xiz]$ and $\cD_K[z]$, 
respectively.  We will also show 
that certain results on $\cA_K[\xiz]$ and 
$\cD_K[z]$ can be extended to the associated 
completions of $\cA_K[\xiz]$ and $\cD_K[z]$.

But, first, let us fix the following notations 
and conventions which 
along with those fixed in the 
previous section will be used 
throughout this paper.

\begin{enumerate}

\item We denote by 
$\langle \cdot, \cdot \rangle$ the standard 
bi-linear paring of $n$-tuples of elements 
of algebras over $K$. For example, 
$\langle \xi, z\rangle\!=\sum_{i=1}^n \xi_i z_i$.

\item We will freely use some standard multi-index notations. 
For instance, let  
$\alpha=(k_1, k_2, ..., k_n)$ and 
$\beta=(m_1, m_2, ..., m_n)$ 
be any $n$-tuples of non-negative 
integers, we have
\begin{align*}
|\alpha|&=\sum_{i=1}^n k_i. \\
\alpha! &=k_1!k_2!\cdots k_n!. \\
\binom{\alpha}{\beta} &=
\begin{cases} \frac{\alpha!}{\beta! (\alpha-\beta)!} 
&\mbox{if  $k_i\ge m_i$ for all $1\le i\le n$}; \\
0, &\mbox{otherwise.}\\
\end{cases}
\end{align*} 

\item For any $k\geq 1$ and
$h(z)=(h_1(z), h_2(z), \cdots, h_k(z)) \in 
\cA_K[[z]]^{\times k}$, 
we define the {\it order}, denoted by $o(h(z))$,  
of $h(z)$ by 
\begin{align*}
o(h(z))\!:=\min_{1\leq i\leq k} o(h_i(z)),  
\end{align*} 
and when $h(z)\in \cA_K[z]^{\times k}$, the {\it degree} 
of $\deg h(z)$ of $h(z)$ by  
\begin{align*}
\deg h(z)\!:=\max_{1\leq i\leq k} \deg h_i(z). 
\end{align*}
For any $h_t(z) \in \cA_{K[t]}[z]^{\times k}
\text{ or } \cA_{K[t]}[z]^{\times k}$ ($k\geq 1$) 
for some formal parameter or variable $t$, 
the notation $o(h_t(z))$ and $\deg h_t(z)$
stand for the {\it order} and 
the {\it degree} of $h_t(z)$ 
with respect to $z$, 
respectively.

\item For any differential operator 
$\phi\in \cD_K[z]$ and polynomial $u(z)\in \cA_K[z]$, 
the notation $\phi\, u(z)$ usually denotes the composition of 
$\phi$ and the multiplication operator by $u(z)$. 
So $\phi \, u(z)$ is still a differential operator 
of $\cA_K[z]$. The polynomial obtained by applying $\phi$ to $u(z)$ 
will be denoted by $\phi(u(z))$ or $\phi\cdot u(z)$. 
\end{enumerate}

Next, we recall some general facts on 
$\bZ$-graded algebras and their 
associated completions.

Let $V$ be a vector space over $K$. 
We say that $V$ is {\it $\bZ$-graded} by 
its subspaces $V_m$ 
$(m\in \bZ)$ if 
$V=\oplus_{m\in \bZ} V_m$, 
i.e., every $v\in V$ can 
be written uniquely as 
$v=\sum_{m\in \bZ} v_m$ 
with $v_m\in V_m$ $(m\in \bZ)$ and all 
but finitely many $v_m=0$. 

For any $\bZ$-graded vector $V$ 
as above, we define the 
{\it associated completion} $\overline{V}$ to be the 
vector space of the elements of the form 
\begin{align}\label{dummy-u}
u=\sum_{m\in \bZ} u_m
\end{align} 
such that  
\begin{enumerate}
  \item[$1)$]  $u_m\in V_m$ for all $m\in \bZ$; 
\item[$2)$] $u_m=0$ when $m$ is large negative enough, 
i.e.,  there exists $N\in \bZ$ (depending on $u$)  
such that $u_m=0$ for all $m\le N$. 
\end{enumerate} 

It is easy to see that the addition and 
the scalar product 
of the vector space $V$ extend to $\overline V$ in 
the obvious way, with which $\overline V$ forms 
a vector space over $K$ and contains $V$ as a 
vector subspace.  

For any $\bZ$-graded vector space $V$ (as above), 
we define the associated {\it grading function}, 
or simply the {\it associated grading} to be the  
function $\eta:\overline{V}\to \bZ \cup \{+\infty\}$ 
such that $\eta(0)=+\infty$ and, for any 
$0\ne u\in \overline{V}$ as in 
Eq.\,(\ref{dummy-u}), 
\begin{align}\label{generic-eta}
\eta(u)=\min \,\{ m \in \bZ\,|\, u_m \ne 0\}.
\end{align}

Note that  for any $m\in \bZ$ and 
$0\ne u\in V_m$, we have 
$\eta(u)=m$. Therefore, we also say  
elements of $V_m$ are 
{\it $\eta$-homogeneous} 
of $\eta$-grading $m$. 
The $\eta$-homogeneous elements $u_m$ 
$(m\in \bZ)$ in Eq.\,(\ref{dummy-u}) 
is called the {\it $\eta$-homogeneous part}  
of the element $u\in \overline{V}$  
of $\eta$-grading $m$. Very often, 
we also say that the vector space 
$\overline{V}$ is the {\it completion} 
of $V$ with respect to the grading 
$\left.\eta\right |_V: 
V \to \bZ \cup\{+\infty\}$.

Now let $\cA$ be any associative 
(but not necessarily commutative) algebra over $K$.
We say  that $\cA$ is {\it $\bZ$-graded} by its subspaces 
$A_m$ $(m\in \bZ)$ if $\cA$ as a vector space 
is $\bZ$-graded by $A_m$ $(m\in \bZ)$ and, 
for any $m, k\in\bZ$, we further have 
\begin{align}
A_m \cdot A_k \subset A_{m+k}.
\end{align}

Let $\bar \cA$ be the associated 
completion of $\cA$ as a $\bZ$-graded 
vector space. Then the product of $\cA$ 
can be extended to $\bar \cA$ in 
the obvious way and, with the extended 
product, $\bar \cA$ forms an associative 
algebra which contains $\cA$ as a 
subalgebra.  

Let $\cA$ be a $\bZ$-graded algebra 
and $\cM$ an $\cA$-module. We say that 
$\cM$ is a {\it $\bZ$-graded} $\cA$-module 
if $\cM$ as a vector space is $\bZ$-graded 
by some of its vector subspaces  
$M_m$ $(m\in \bZ)$ such that   
for all $m, k\in\bZ$, we have 
\begin{align}
A_m \cdot M_k\subset M_{m+k}.
\end{align}

Furthermore, let $\overline{M}$ be the 
associated completion of $\cM$. 
Then the $\cA$-module structure 
on $\cM$ can be extended naturally to 
an $\cA$-module structure on 
$\overline{M}$. 
This $\cA$-module structure on 
$\overline{M}$ can be further 
extended in the obvious way to an $\bar \cA$-module 
structure. 

\begin{exam}\label{E2.1}
$a)$ The polynomial algebra $\cA_K[z]$  
is $\bZ$-graded by the degree of polynomials in 
$\cA_K[z]$. More precisely, for any $m\in \bZ$, set 
$A_m[z]=0$ if $m<0$, and $A_m[z]$ the subspace 
of homogeneous polynomials in $z$ of degree $m$. 
Then it is easy to see that $\cA_K[z]$ is $\bZ$-graded 
by the subspaces $A_m[z]$ $(m\in \bZ)$ and 
the associated completion 
$\bar \cA_K[z]$ of $\cA_K[z]$ is the formal power series 
algebra $\cA_K[[z]]$ in $z$ over $K$.

$b)$ Let $\cA_K[z, z^{-1}]$ denote  
the algebra of Laurent polynomials 
in $z$ over $K$. 
Then $\cA_K[z, z^{-1}]$ is also $\bZ$-graded by 
the $($generalized$)$ degree of Laurent polynomial 
in $z$, i.e.,  counting the degree of $z_i^{-1}$ by $-1$ 
for each $\lin$. The associated completion of 
$\cA_K[z, z^{-1}]$ with respect to this $\bZ$-grading 
is the algebra of Laurent power 
series in $z$ over $K$.
\end{exam} 	

Besides the one in the example $a)$ above, in this paper 
we also need to work over 
the following two $\bZ$-graded algebras. 
 
The first one is the polynomial algebra $\cA_K[\xiz]$ over $K$ 
in the $2n$ variables $\xi_i$ and $z_i$ $(\lin)$ fixed 
at the beginning of this section. 
For any $m\in \bZ$, we set 
\begin{align}\label{eta-hgs}
A_m[\xiz]\!:=\mbox{Span}_K \{\xi^\alpha z^\beta\,|\,
\alpha, \beta \in \bN^n; |\beta|-|\alpha|=m\}.
\end{align}

Then  it is easy to see that 
$\cA_K[\xiz]$ is $\bZ$-graded by 
the subspaces $A_m[\xiz]$ 
$(m\in \bZ)$ and its associated completion, 
denoted by $\bar \cA_K[\xiz]$, is a proper 
subalgebra of the formal power series $\cA_K[[\xiz]]$.
Note that the associated completion 
$\bar \cA_K[\xiz]$ contains the formal 
power series algebra $\cA_K[[z]]$ in 
$z$ as a subalgebra. 
But it does not contain the formal 
power series algebra 
$\cA_K[[\xi]]$ in $\xi$.
 
Let $\eta: \bar \cA_K[\xiz]\to \bZ\cup\{+\infty\}$
be the associated $\bZ$-grading of 
$\cA_K[\xiz]$. Then  for any $\alpha, \beta\in \bN^n$, 
we have
\begin{align}\label{Def-eta-grading}
\eta(\xi^\alpha z^\beta)=|\beta|-|\alpha|.
\end{align}
More generally, for any $f(\xiz)\in \bar \cA[\xiz]$, 
\begin{align}\label{Def-eta-grading-2}
\eta(f(\xiz))=\min\{ |\beta|-|\alpha|\,|\, 
\alpha, \beta\in \bN^n; [\xi^\alpha z^\beta]f(\xiz)\ne 0\}, 
\end{align}
where $[\xi^\alpha z^\beta]f(\xiz)$ denotes the 
coefficient of the monomial $\xi^\alpha z^\beta$
in $f(\xiz)$.

The second $\bZ$-graded algebra that we will need 
later is the Weyl algebra $\cD_K[z]$ of 
the polynomial algebra $\cA_K[z]$. 

Note first that it is well-known (e.g., by using 
Proposition 2.2 (p.\,4) in \cite{B} 
or Theorem 3.1 (p.\,58) in \cite{C}) 
that the differential operators 
$z^\beta\p^\alpha$  $(\alpha, \beta\in \bN^n)$
form a linear basis of $\cD_K[z]$. 
So do the differential operators $\p^\alpha z^\beta$  
$(\alpha, \beta\in \bN^n)$.

For any $m\in \bZ$, we set 
\begin{align}\label{nu-hgs}
D_m[z]\!:=\mbox{Span}_K \{z^\beta \p^\alpha 
\,|\,\alpha, \beta \in \bN^n; |\beta|-|\alpha|=m\}.
\end{align}

Then  it is easy to see that the Weyl algebra $\cD_K[z]$ 
is $\bZ$-graded by the subspaces $D_m[z]$ 
$(m\in \bZ)$. Actually, this grading is the same 
as the grading in the Kashiwara-Malgrange $V$-filtration 
along the origin (see \cite{B2}, \cite{D}). 

We denote by $\overline{\cD}_K[z]$ 
the associated completion of $\cD_K[z]$ and  
$\nu: \overline{\cD}_K[z] \to \bZ\cup\{+\infty\}$
the associated $\bZ$-grading of 
$\cD_K[z]$. 
Then  for any $\alpha, \beta\in \bN^n$, 
we have 
\begin{align}\label{nu-grading-1}
\nu(z^\beta\p^\alpha)=|\beta|-|\alpha|.
\end{align}
 
Note also that for any $\alpha, \beta\in \bN^n$, 
we also have 
\begin{align}\label{nu-grading-2}
\nu(\p^\alpha z^\beta)=|\beta|-|\alpha|.
\end{align}

This is because by the Leibniz rule, 
we have  
$$
\p^\alpha z^\beta=
\sum_{\gamma\in \bN^n} 
\binom{\alpha}{\gamma} 
(\p^\gamma(z^\beta))\p^{\alpha-\gamma},
$$
and each differential operator 
$(\p^\gamma(z^\beta))\p^{\alpha-\gamma}$ with a 
non-zero coefficient in the sum above 
has $\nu$-grading $|\beta|-|\alpha|$.

With the $\bZ$-grading of the polynomial 
algebra $\cA_K[z]$ as in Example \ref{E2.1}, $a)$
and the $\bZ$-grading of the Weyl 
algebra $\cD_K[z]$ defined above, it is easy to see 
that $\cA_K[z]$ is a $\bZ$-graded $\cD_K[z]$-module. 
Therefore, the associated 
completion $\bar \cA_K[z](=\cA_K[[z]])$ of 
$\cA_K[z]$ becomes a module 
of the associated completion $\overline \cD_K[z]$
of the Weyl algebra $\cD_K[z]$. 
In other words, we have a well-defined action 
of $\overline \cD_K[z]$ on  
the formal power series 
algebra $\cA_K[[z]]$.

Next, we consider the {\it right} and 
{\it  left total symbols} of the differential operators 
in the associated completion $\overline\cD_K[z]$ 
of the Weyl algebra $\cD_K[z]$. 

First let us recall the {\it right} and 
{\it  left total symbols} 
of the differential operators in $\cD_K[z]$, i.e.,  
the differential operators of the polynomial algebra 
$\cA_K[z]$.   

For any $\phi \in \cD_K[z]$,  
we can write $\phi$ uniquely as the following two finite sums:
\begin{align}\label{generic-diff}
\phi=\sum_{\alpha \in \bN^n} a_\alpha(z) \p^\alpha
=\sum_{\beta \in \bN^n} \p^\beta b_\beta(z)
\end{align}  
where $a_\alpha(z), b_\beta(z) \in \cA_K[z]$ but  
denote the multiplication operators by 
$a_\alpha(z)$ and $b_\beta(z)$, respectively. 

For the differential operator $\phi\in \cD_K[z]$ 
in Eq.\,(\ref{generic-diff}), 
its {\it right} and {\it left total symbols} are defined 
to be the polynomials $\sum_{\alpha\in \bN^n} a_\alpha(z) 
\xi^\alpha\in \cA_K[\xiz]$ and 
$\sum_{\beta \in \bN^n} b_\beta(z) 
\xi^\beta\in \cA_K[\xiz]$, respectively.
We denote by $\cR: \cD[z] \to \cA_K[\xiz]$
(resp., $\cL: \cD[z] \to \cA_K[\xiz]$) the linear map 
which maps any $\phi\in \cD[z]$ to its right 
(resp., left) total symbol.

Furthermore, from the definitions of the 
$\bZ$-gradings of $\cA_K[\xiz]$ and $\cD_K[z]$, 
and also Eq.\,(\ref{nu-grading-2}),   
it is easy to see that both $\cR$ and $\cL$ are 
$\bZ$-grading preserving, i.e.,  for any $m\in \bZ$, 
we have
\begin{align}
\cR (D_m[z]) &\subset A_m[\xiz]. \\ 
\cL (D_m[z]) &\subset A_m[\xiz]. 
\end{align}

Since $\cR$ and $\cL$ 
are isomorphisms, for each $m\in \bZ$  
the restrictions of $\cR$ and $\cL$ on 
$D_m[z]$ are also isomorphisms from $D_m[z]$ 
to $A_m[\xiz]$. Hence, the  
following lemma holds. 

\begin{lemma}\label{RL-lemma}
$\cR$ and $\cL$ extend 
to $\bZ$-grading preserving linear isomorphisms 
from the completion $\overline\cD_K[z]$ of $\cD_K[z]$ 
to the completion $\bar\cA_K[\xiz]$ 
of $\cA_K[\xiz]$. 
\end{lemma}

We will still denote 
by $\cR$ and $\cL$ the extended isomorphisms 
of $\cR$ and $\cL$, respectively. For any 
$\phi\in \overline\cD_K[z]$, the images 
$\cR(\phi)$ and $\cL(\phi)$ will also 
be called respectively the {\it right} and 
{\it left total symbols} of $\phi$.

Next we introduce the following linear 
automorphism $\Phi$ of $\cA_K[\xiz]$.
Set 
\begin{align}
\Lambda\!:&=\sum_{i=1}^n \p_{\xi_i}\p_{z_i}, \label{Def-Lambda} \\
\Phi\!:&=e^\Lambda=\sum_{m \ge 0} \frac{\Lambda^m}{m!}.
\label{Def-Phi}
\end{align}

First, note that $\Phi$ is a well-defined linear map 
from $\cA_K[\xiz]$ to $\cA_K[\xiz]$. 
This is because  $\Lambda$ is locally nilpotent, 
i.e.,  for any $f(\xiz)\in \cA_K[\xiz]$, $\Lambda^m f=0$ 
when $m\gg 0$. Also, $\Phi$ is invertible with 
the inverse map given by $e^{-\Lambda}$.

Second, it is easy to see that 
$\Lambda$ is a linear endomorphism 
of $\cA_K[\xiz]$ and preserves the 
$\eta$-grading of 
$\cA_K[\xiz]$ given in 
Eqs.\,(\ref{Def-eta-grading}) and 
(\ref{Def-eta-grading-2}). 
Hence, so does $\Phi$.
Therefore, $\Phi$ extends to 
an automorphism, still denoted by 
$\Phi$, of the associated completion 
$\bar\cA_K[\xiz]$ of $\cA_K[\xiz]$ .  

\begin{propo}\label{Phi-Trans}
For the extended isomorphisms $\Phi$, $\cR$ and $\cL$ 
introduced above, 
we have $\Phi=\cR \circ \cL^{-1}$ 
$($as linear automorphisms of $\bar\cA_K[\xiz]$$)$.   
\end{propo}

\pf Note first that  by Lemma \ref{RL-lemma} and 
the discussion before the proposition, 
we know that both $\Phi$ and 
$\cR \circ \cL^{-1}$ are $\eta$-grading 
preserving automorphisms 
of $\bar\cA_K[\xiz]$.

Second, by Theorem $3.1$ in \cite{T-Deform}, 
we know that the restrictions of 
$\Phi$ and $\cR \circ \cL^{-1}$ on 
the subalgebra $\cA_K[\xiz]\subset \bar\cA_K[\xiz]$
coincide. In particular, their restrictions on 
the $\eta$-homogeneous subspaces 
$A_m[\xiz]$ $(m\in \bZ)$ in Eq.(\ref{eta-hgs}) 
of $\cA_K[\xiz]$ coincide.  
Hence, as linear automorphisms of 
$\bar\cA_K[\xiz]$, $\Phi$ and 
$\cR \circ \cL^{-1}$ must also 
be same.   
\epfv

Note that  even though the base field in 
\cite{T-Deform} is the complex number field $\bC$, 
the proof of Theorem $3.1$ in 
\cite{T-Deform} works equally well for all the fields 
of characteristic zero. 

Finally, let us give the following lemma that will be needed later.

Let $\tau:\cD_K[z]\to \cD_K[z]$ be the linear map 
defined by setting 
\begin{align}\label{Def-tau}
\tau( h(z) \p^\alpha )=(-1)^{|\alpha|} 
\p^\alpha h(z)
\end{align} 
for each $h(z)\in \cA_K[z]$ and 
$\alpha\in \bN^n$. 

Then it is well-known (e.g., see $\S16.2$ in \cite{C}) that $\tau$ is an anti-involution of the Weyl algebra 
$\cD_K[z]$, i.e.,    
$\tau^2=id$ and, for any 
$\phi, \psi\in \cD_K[z]$, 
we have
\begin{align}\label{anti-tau}
\tau(\phi\psi) = \tau(\psi) \tau(\phi).
\end{align}

Note that  by Eqs.\,(\ref{nu-grading-1}) and 
(\ref{nu-grading-2}), we see that 
$\tau$ preserves the $\nu$-grading of 
$\cD_K[z]$. Hence, $\tau$ extends to 
an automorphism, 
still denoted by $\tau$,  
of the associated completion 
$\overline \cD_K[z]$ 
of the weyl algebra $\cD_K[z]$. 
Furthermore, the following lemma can also 
be easily checked. 

\begin{lemma}\label{tau-lemma}
The extended $\tau$ is also an anti-involution of 
the completion $\overline{\cD}_K[z]$ of 
the Weyl algebra $\cD_K[z]$.
\end{lemma}

\renewcommand{\theequation}{\thesection.\arabic{equation}}
\renewcommand{\therema}{\thesection.\arabic{rema}}
\setcounter{equation}{0}
\setcounter{rema}{0}

\section{\bf The Abhyankar-Gurjar Inversion Formula Revisited}\label{S3}

In this section, we first give a new proof and 
then a new formulation 
(see  Theorem \ref{AG-in-Phi}) for 
the Abhyankar-Gurjar inversion formula (see  Theorem \ref{AGF}). 
Some consequences of the reformulated Abhyankar-Gurjar inversion formula will also be derived. But, first, let us fix the following notation.

Let $F(z)=(F_1, F_2, ... , F_n)$ be a homomorphism of  
the formal power series algebra $\cA_K[[z]]$ such that $F(0)=0$. 
Up to some conjugations by linear automorphisms, 
we may and will always assume $F$ has the form $F(z)=z-H(z)$ and 
$H(z)=(H_1, H_2, ... , H_n)\in\cA_K[[z]]^{\times n}$ with the order $o(H) \ge 2$. 

With the assumption above on $F$, we have $(\cJ F)(0)=1$. Therefore, $F$ has a formal inverse map $G(z)=F^{-1}(z)$ 
which can be written as $G(z)=z+N(z)$ for some 
$N(z)=(N_1, N_2, ..., N_n) \in \cA_K[[z]]^{\times n}$ 
with  $o(N) \ge 2$.

With the setup above, the Abhyankar-Gurjar  inversion 
formula can be stated as follows.

\begin{theo} \label{AGF} $(${\bf Abhyankar-Gurjar}$)$
For any $u(z)\in \cA_K[[z]]$, we have 
\begin{align}\label{AGF-e1}
u(G(z))= \sum_{\alpha \in \bN^n} \frac 1{\alpha!}\, \p^\alpha 
\left( u H^\alpha \cJ F \right)
=\Big(\sum_{\alpha \in \bN^n} \frac 1{\alpha!}\, \p^\alpha 
  H^\alpha \cJ F \Big) \cdot u(z).
\end{align}

In particular, the formal inverse map $G$ of $F$ is 
given by 
\begin{align}\label{AGF-e2}
G(z)= \sum_{\alpha \in \bN^n} 
\frac 1{\alpha!} \p^\alpha \, 
\left( z H^\alpha \cJ F  \right).
\end{align}
\end{theo}

\pf First, by Proposition $2.1$ and Lemma $2.4$ in \cite{Exp}, 
we know that there exists a unique 
$a(z)=(a_1(z), a_2(z), ..., a_n(z))\in \cA_K[[z]]^{\times n}$ 
with $o(a(z))\ge 2$, such that  for any $u(z)\in \cA_K[[z]]$, 
\begin{align}
e^{a(z)\p}(u(z))=u(F), \label{AGF-pe1} \\
e^{-a(z)\p}(u(z))=u(G), \label{AGF-pe1b}
\end{align}
where $a(z)\p\!:=\sum_{i=1}^n a_i(z)\p_i$.

On the other hand, by the Taylor series 
expansion of $u(F)=u(z-H)$ at $z$, we have
\begin{align}\label{AGF-pe2} 
u(F)&=\sum_{\alpha\in \bN^n} 
\frac {(-1)^{|\alpha|}}{\alpha!}  
H^\alpha (\p^\alpha u)(z)= \left( \sum_{\alpha\in \bN^n} 
\frac {(-1)^{|\alpha|}}{\alpha!}  
H^\alpha \p^\alpha\right) \cdot u(z). 
\end{align}

Note that  by the facts that $o(a)\ge 2$ and $o(H)\ge 2$,
it is easy to see that both the differential operators
$e^{a(z)\p}$ and $\sum_{\alpha\in \bN^n} 
\frac {(-1)^{|\alpha|}}{\alpha!}  
H^\alpha(z) \p^\alpha$ lie in the completion 
$\overline\cD_K[z]$ of $\cD_K[z]$ constructed 
in the previous section.
Furthermore, since Eqs.\,(\ref{AGF-pe1}) 
and (\ref{AGF-pe2})
hold for all $u(z)\in\cA_K[[z]]$, 
hence we have 
\begin{align}\label{AGF-pe3} 
e^{a(z)\p}&=\sum_{\alpha\in \bN^n} 
\frac {(-1)^{|\alpha|}}{\alpha!}  
H^\alpha \p^\alpha. 
\end{align}

Next, we consider the images of 
the differential operators in the 
equation above under the anti-involution 
$\tau$ of $\overline\cD_K[z]$ 
in Lemma \ref{tau-lemma}.

First, by Eq.\,(\ref{Def-tau}) we have 
\begin{align}\label{AGF-pe4} 
\tau(a(z)\p)&=-\sum_{i=1}^n \p_i a_i(z)= 
-\sum_{i=1}^n (a_i(z)\p_i +\p_i ( a_i(z) ) \\
&=-a(z)\p-\nabla a(z), \nno
\end{align}  
where $\nabla a(z)\!:=\sum_{i=1}^n (\p_i a_i)(z)$, i.e.,  the {\it gradient} of the $n$-tuple $a(z)$, but it here denotes the multiplication operator by $\nabla a$.

Second, by Lemma \ref{tau-lemma}, $\tau$ is also an anti-involution of 
$\overline\cD_K[z]$, whence we have 
\begin{align}\label{AGF-pe5} 
\tau(e^{a(z)\p})=e^{\tau(a(z)\p)}= 
e^{-a(z)\p-\nabla a}.
\end{align}
  
Now, apply $\tau$ to Eq.\,(\ref{AGF-pe3}),   
by Eqs.\,(\ref{AGF-pe5}) and (\ref{Def-tau}), we have 
\begin{align}\label{AGF-pe6} 
e^{-a(z)\p-\nabla a(z)}=
\sum_{\alpha\in \bN^n} \frac 1{\alpha!} 
\p^\alpha H^\alpha. 
\end{align}

Note that  for any $u(z)\in \cA_K[[z]]$, 
by Theorem $2.8$, $b)$ in \cite{Exp} with 
$F(t, z)$ there replaced by $G(z)$, or 
equivalently, by letting $t=-1$,   
we have 
\begin{align}\label{AGF-pe7} 
e^{-a(z)\p-\nabla a(z)} (u(z))= \cJ G (z) u(G(z)).
\end{align}

Combining Eqs.\,(\ref{AGF-pe6}) and (\ref{AGF-pe7}), 
we have    
\begin{align}\label{AGF-pe8}
\sum_{\alpha\in \bN^n} \frac 1{\alpha!} \, 
\p^\alpha (H^\alpha u(z))= \cJ G (z) u(G(z))
\end{align}
for all $u(z)\in \cA_K[[z]]$.

By replacing $u(z)$ in Eq.\,(\ref{AGF-pe8}) by 
$u(z)\cJ F(z)$, and by the identity 
$\cJ F(G)\cJ G (z)=\cJ(F(G))=\cJ z=1$, 
we have  
\begin{align*}
\sum_{\alpha\in \bN^n} \frac 1{\alpha!} 
\p^\alpha (H^\alpha \cJ F(z) u(z))=\cJ G(z) 
\cJ F(G) u(G(z))
=u(G(z)), 
\end{align*}
which is Eq.\,(\ref{AGF-e1}). 
\epfv

Next we give a reformulation of the Abhyankar-Gurjar 
inversion formula in terms of certain elements of $\bar\cA_K[\xiz]$ 
and the differential operator $\Phi =e^\Lambda$ 
defined in Eqs.\,(\ref{Def-Lambda}) and (\ref{Def-Phi}).

\begin{theo}\label{AG-in-Phi}
Let $H$ and $N$ be as above. Then  
for any $q(z)\in \cA_K[[z]]$, 
we have
\begin{align}
\Phi \left( q(z) (\cJ F)(z) e^{\langle \xi, H(z)\rangle} \right)
&=q(G) \, e^{\langle \xi, N(z)\rangle}. \label{AG-in-Phi-e1} 
\end{align}
\end{theo}

Note that  the formal power series 
$q(z) (\cJ F)(z) e^{\langle \xi, H(z)\rangle}$ on the left hand side of Eq.\,(\ref{AG-in-Phi-e1}) lies in $\bar\cA_K[\xiz]$ due to the assumption $o(H)\ge 2$. So, by the discussion in the previous section, the left hand side of 
Eq.\,(\ref{AG-in-Phi-e1}) makes sense 
and gives an element of $\bar\cA_K[\xiz]$. 
Similarly, due to the fact $o(N)\ge 2$, the element 
$q(G)e^{\langle \xi, N(z)\rangle}$ on right hand side of 
Eq.\,(\ref{AG-in-Phi-e1}) also lies 
in $\bar\cA_K[\xiz]$. \\

{\underline{\it Proof of Theorem \ref{AG-in-Phi}:}\, 
First, for any $u(z)\in \cA_K[[z]]$, 
by Eq.\,(\ref{AGF-e1}) with $u(z)$ replaced 
by $q(z)u(z)$, we have 
\begin{align}
q(G)u(G)=\sum_{\alpha \in \bN^n} \frac 1{\alpha!} \p^\alpha 
\big( H^\alpha(z) \cJ F(z) q(z) u(z)\big). 
\end{align}

Second, by the Taylor series of $u(G)=u(z+N)$ at $z$, we have 
\begin{align}
q(G)u(G)&=q(G)\sum_{\alpha \in \bN^n}\frac 1{\alpha!}N^\alpha(z) 
\p^\alpha (u(z)) \\
&= \sum_{\alpha \in \bN^n}\frac 1{\alpha!} q(G) N^\alpha(z) 
\p^\alpha (u(z)). \nno 
\end{align}

Therefore, as differential operators in 
$\overline\cD_K[z]$, we have 
\begin{align}
\sum_{\alpha \in \bN^n} \frac 1{\alpha!} \, 
\p^\alpha H^\alpha(z) \cJ F(z) q(z) = 
\sum_{\alpha \in \bN^n}\frac 1{\alpha!} \, 
q(G) N^\alpha(z) \p^\alpha. 
\end{align}

In terms of the left and right total symbol 
maps $\cL$ and $\cR$ from $\overline\cD_K[z]$ 
to $\bar\cA_K[\xiz]$ (see lemma \ref{RL-lemma}), 
the equation above can be written further 
as follows:
\allowdisplaybreaks{
\begin{align}
\cL^{-1} \left(\sum_{\alpha \in \bN^n} \frac 1{\alpha!} \, 
\xi^\alpha H^\alpha(z) (\cJ F)(z) q(z) \right)& = 
\cR^{-1}\left(
\sum_{\alpha \in \bN^n}\frac 1{\alpha!} \, 
q(G) N^\alpha(z) \xi^\alpha \right). \nno 
\end{align}
Take sum over $\alpha\in \bN^n$ and 
then apply $\cR$ to the both sides, we get 
\begin{align}
\cL^{-1} \left( q(z) (\cJ F)(z) 
e^{\langle \xi, H(z)\rangle} \right)& = 
\cR^{-1}\left(  
q(G) e^{\langle \xi,N(z)\rangle} \right). \nno \\
(\cR\circ \cL^{-1}) \left( q(z) (\cJ F)(z) 
e^{\langle \xi, H(z)\rangle} \right)& =   
q(G) e^{\langle \xi,N(z)\rangle}.\label{T3.2-lasteq}
\end{align} }
But, by Proposition \ref{Phi-Trans} 
we also have $\cR\circ \cL^{-1}=\Phi$. 
Hence, Eq.\,(\ref{AG-in-Phi-e1}) follows.
\epfv

Note that  by letting $q(z)=1$ in Eq.\,(\ref{AG-in-Phi-e1}), 
we immediately have the following corollary.

\begin{corol}\label{Lambda-Inversion}
Let $F$ and $H$ be fixed as before. Then  
\begin{align}\label{Lambda-Inversion-e1}
\Phi \left( \cJ F (z) 
e^{\langle \xi, H(z)\rangle}\right)
=e^{\langle \xi, N(z)\rangle}.
\end{align} 

In particular, when $\cJ F\equiv 1$, we have
\begin{align}\label{Lambda-Inversion-e2}
\Phi \left( e^{\langle \xi, H(z)\rangle}\right)
&=e^{\langle \xi, N(z)\rangle}. 
\end{align} 
\end{corol}

\begin{rmk}
By going backward the proof of 
Theorem \ref{AG-in-Phi} with $q(z)=1$,
it is easy to see that Eq.\,$(\ref{AGF-e1})$ 
also follows from Eq.\,$(\ref{Lambda-Inversion-e1})$.
Therefore,  Eq.\,$(\ref{Lambda-Inversion-e1})$ or 
Eq.\,$(\ref{AG-in-Phi-e1})$ can be viewed as 
a reformulation of the Abhyankar-Gurjar inversion 
formula in Eq.\,$(\ref{AGF-e1})$.
\end{rmk}

By using Theorem \ref{AG-in-Phi}, we can also derive 
the following formula, which will play a crucial 
role in next section. 

\begin{theo}\label{Lambda-Inversion-2}
Let $P=\langle \xi, H\rangle$. Then  
for any $k\ge 0$ and $q(z)\in \cA_K[z]$, we have 
\begin{align}\label{Lambda-Inversion-e3'}
q(G(z))\langle \xi, N(z)\rangle^k=
k!\sum_{m\ge 0} \frac{\Lambda^m \left( P^{m+k} q(z) \,\cJ F \right)}{m!(m+k)!}.
\end{align}

In particular, by choosing $q(z)=1$,  we get 
\begin{align}\label{Lambda-Inversion-e3}
\langle \xi, N(z)\rangle^k=
k!\sum_{m\ge 0} \frac{\Lambda^m \left( P^{m+k}\cJ F\right)}{m!(m+k)!}.
\end{align}
\end{theo}

\pf We view the formal power series of $\xi$ and $z$ on the both sides of Eq.\,(\ref{AG-in-Phi-e1}) as formal power 
series in $\xi$ with coefficients in $\cA_K[[z]]$. Then  
by comparing the homogeneous parts in $\xi$ of 
degree $k$ of both sides of Eq.\,(\ref{AG-in-Phi-e1}) 
and noting that $\Phi=e^\Lambda$, 
we have 
\begin{align*}
\sum_{m\ge 0} \frac{\Lambda^m \left( P^{m+k}q(z)\,\cJ F\right)}{m!(m+k)!}
=\frac 1{k!}\, q(G(z))\, \langle \xi, N(z)\rangle^k,
\end{align*}
which gives us  
Eq.\,(\ref{Lambda-Inversion-e3'}).
\epfv

Note also that by letting $k=0$ in 
Eq.\,(\ref{Lambda-Inversion-e3'}), we immediately 
have the following inversion formula. 

\begin{corol}\label{New-Inv}
Let $P$ be as in Theorem \ref{Lambda-Inversion-2}. 
Then for any $q(z)\in \cA_K[[z]]$, we have 
\begin{align}\label{New-Inv-e1}
q(G(z))=\sum_{m\ge 0} \frac{\Lambda^m \left(P^{m}q(z)\, \cJ F \right)}{(m!)^2}.
\end{align}

In particular, by choosing $q(z)=z$, we get 
\begin{align}\label{New-Inv-e2}
G(z)=\sum_{m\ge 0} \frac{\Lambda^m \left( z P^{m}\cJ F \right)}{(m!)^2}.
\end{align}
\end{corol}

\renewcommand{\theequation}{\thesection.\arabic{equation}}
\renewcommand{\therema}{\thesection.\arabic{rema}}
\setcounter{equation}{0}
\setcounter{rema}{0}

\section{\bf A New Proof For the Equivalence of the Jacobian Conjecture and the Vanishing Conjecture}\label{S4}

In this section, we show the equivalence 
of the  {\it Jacobian conjecture}  with a special case of  
the  {\it vanishing conjecture}, 
Conjecture \ref{VC}. But,  we first need to recall  
the following deformation of formal maps.

Let $F=z-H$ and $G=z+N$ as fixed in the previous section. 
Let $t$ be a formal parameter which commutes with 
$z_i$ $(1\leq i\leq n)$. We set 
$F_t(z)\!:=z-tH(z)$. Since $F_{t=1}(z)=F(z)$, 
$F_t(z)$ can be viewed as a 
deformation of the formal map $F(z)$.

We may view $F_t(z)$ as a formal map of 
$\cA_{K(t)}[[z]]$, i.e.,  a formal map 
over the rational function field $K(t)$. 
Hence, with base 
field $K$ replaced by $K(t)$, all 
results derived in the previous 
sections apply to 
the formal map $F_t(z)$. 

Denote by $G_t(z)$ the formal inverse  
of $F_t(z)$. It is well known that    
$G_t(z)$ can actually be written as 
$G_t(z)=z+tN_t(z)$ for some 
$N_t(z)\in \cA_{K[t]}[[z, t]]^{\times n}$  
with $\mbox{o} (N_t(z))\geq 2$, i.e.,  all 
the coefficients of $N_t(z)$ 
are actually polynomials of $t$ over $K$.
For example, this can be easily seen  
by applying the Abhyankar-Gurjar 
inversion formula in Eq.\,(\ref{AGF-e2}) 
to the formal map $F_t(z)$. 

Another remark is that  since $F_{t=1}(z)=F(z)$, 
by the uniqueness of inverse maps and the fact 
$N_t(z)\in \cA_{K[t]}[[z]]^{\times n}$ mentioned 
above, it is easy to see that  we also have   
$G_{t=1}(z)=G(z)$.

Next, let us start with the following lemma.

\begin{lemma}\label{H-nil-lemma}
$JH$ is nilpotent iff $\cJ F_t\equiv 1$.
\end{lemma}

\pf First, since $\cA_K[[z]]$ is an integral domain over $K$, as a matrix with entries in $\cA_K[[z]]$, $JH$ is nilpotent iff 
its characteristic polynomial $\det(\lambda I-JH)=\lambda^n$.

Second, $\cJ F_t=\det(I-t JH)=t^n\det(t^{-1}I-JH)$. 
Therefore, $\cJ F_t\equiv 1$ iff $\det(t^{-1}I-JH)=t^{-n}$ 
iff the characteristic polynomial 
of $JH$ is equals to $\lambda^n$. 
Hence the lemma follows.
\epfv

\begin{propo}\label{Lambda-Nil}
Let $P(\xiz)=\langle \xi, H(z)\rangle$. Then   
$JH$ is nilpotent iff $\Lambda^m(P^m)=0$ 
for  all $m\ge 1$.
\end{propo}

\pf 
Applying Eq.\,(\ref{AG-in-Phi-e1}) to the formal map $F_t$ 
with $q(z)=1/\cJ F_t(z)$, we have 
\begin{align}
 e^{\Lambda}  (e^{t \langle \xi, H(z)\rangle})
= \cJ G_t(z) e^{t\langle \xi, N_t(z)\rangle}.
\end{align}
Applying the change of variables $z\to z$ and 
$\xi\to t^{-1}\xi$ in the equation above, and  
noting that $\Lambda\to t\Lambda$, 
we get 
\begin{align}
 e^{t\Lambda} (e^{\langle \xi, H(z)\rangle})
= \cJ G_t(z) e^{\langle \xi, N_t(z)\rangle}.
\end{align}
By looking at the homogeneous part in $\xi$ 
of degree zero of the equation above,
we get 
\begin{align}\label{Lambda-Nil-pe3}
\sum_{m\ge 0} \frac{t^m \Lambda^m (P^{m})}{(m!)^2}
=\cJ G_t(z).
\end{align}

From the equation above, it is easy to see 
that $\cJ G_t(z)\equiv 1$ iff 
$\Lambda^m (P^{m})=0$ for  all  $m\ge 1$.

On the other hand, since 
$\cJ F_t(z)=1/(\cJ G_t)(F_t)$, 
we also have that  $\cJ F_t(z)\equiv 1$ 
iff $\cJ G_t(z)\equiv 1$.

Hence, we have  $\cJ F_t(z)\equiv 1$ iff 
$\Lambda^m (P^{m})=0$ for  all  $m\ge 1$.
Then the proposition follows immediately 
from this fact and 
Lemma \ref{H-nil-lemma}
\epfv
 
When $F=z-H$ is a polynomial map, the $(\Leftarrow)$ part of Proposition \ref{Lambda-Nil} can actually be improved as follows. 

\begin{propo}\label{Lambda-Nil-2}
Let $F=z-H$ and $P$ as in 
Proposition \ref{Lambda-Nil} but 
with $H(z)\in \cA_K[z]^{\times n}$. 
Assume that $\Lambda^m (P^m)=0$ when $m\gg 0$. 
Then  $JH$ is nilpotent.
\end{propo}

\pf First, by Eq.\,(\ref{Lambda-Nil-pe3}) and the assumption 
of the proposition, we see that $\cJ G_t(z)$ is actually 
a polynomial in $z$ and $t$, i.e.,  $\cJ G_t(z)\in \cA_K[z, t]$. 

Second, since $H(z)\in \cA_K[z]^{\times n}$ and $F_t=z-tH$, 
both $\cJ F_t=\det(I-tJH)$ and the composition $(\cJ G_t)(F_t)$ are 
also in $\cA_K[z, t]$. 

Since the polynomial algebra $\cA_K[z, t]$ is an integral domain, 
from the identity $(\cJ G_t)(F_t)\,\cJ F_t=1$ we get 
$\cJ F_t=1$. Then by Lemma \ref{H-nil-lemma},   
$JH$ is nilpotent. 
\epfv

Note that  combining with Proposition \ref{Lambda-Nil}, 
Proposition \ref{Lambda-Nil-2} 
gives a positive answer for a special case of the 
following conjecture.

\begin{conj}\label{Stabe-Nil}
For any $P(\xiz)\in \cA_K[\xiz]$ with 
$\Lambda^m (P^m)=0$ when $m\gg 0$, we have 
$\Lambda^m (P^m)=0$ for  all  $m\ge 1$.
\end{conj}

Note that  up to changes of variables, the conjecture above 
is actually equivalent to Conjecture $4.4$ in \cite{HNP},   
whose general case is still open.

\begin{propo}\label{t-Lambda-JC}
Let $F=z-H$ and $P$ as above. Assume further that 
$F$ is a polynomial map with $JH$ nilpotent. 
Then  the following two statements are equivalent.

$(a)$ $G_t(z)$ is a polynomial map over $K[t]$.

$(b)$ $\Lambda^m(P^{m+1})=0$ when $m\gg 0$. 
\end{propo}

\pf Note that since $JH$ is nilpotent, 
$\cJ F_t(z)=\det(I-tJH)\equiv 1$.  
Applying Eq.\,(\ref{Lambda-Inversion-e3}) 
with $k=1$ to $F_t$, we have
\begin{align*}
\langle \xi, tN_t(z)\rangle
= \sum_{m\ge 0} 
\frac{\Lambda^m (\langle \xi, tH\rangle^{m+1})}{m!(m+1)!}
=\sum_{m\ge 0} \frac{t^{m+1}\Lambda^m (P^{m+1})}{m!(m+1)!}.
\end{align*}

Hence, we have   
\begin{align}\label{t-Lambda-JC-pe1}
\langle \xi, N_t(z)\rangle
=\sum_{m\ge 0} \frac{t^m \Lambda^m (P^{m+1})}{m!(m+1)!}.
\end{align}

From the equation above, it is easy to see that,   
$N_t(z)$ is a $n$-tuple of polynomials in $t$ 
with coefficients in $\cA_K[z]$ iff  
$\Lambda^m (P^{m+1})=0$ when $m\gg 0$.

On the other hand, as pointed out at the beginning of this section, 
$N_t(z) \in \cA_{K[t]}[[z]]^{\times n}$. Therefore,
$N_t(z)$ is a $n$-tuple of polynomials in $t$ 
with coefficients in $\cA_K[z]$ iff 
$N_t(z)$ is a $n$-tuple of polynomials in $z$ 
with coefficients in $K[t]=\cA_K[t]$.

Combining the two equivalences above, we get the 
equivalence of the statements $(a)$ and $(b)$ 
in the proposition.
\epfv

Now we are ready to formulate and prove the main result of this section.

\begin{theo}\label{JC=VC}
The following two statements are equivalent to each other.

$(a)$ The  {\it Jacobian conjecture}  holds for all $n\ge 1$.

$(b)$ For any $n\ge 1$ and  
$H(z)\in \cA_K[z]^{\times n}$ with $o(H)\ge 2$, 
Conjecture \ref{VC} holds for the quadratic 
differential operator  
$\Lambda=\sum_{i=1}^n \p_{\xi_i}\p_{z_i}$ 
and the polynomial 
$P(\xiz)\!:=\langle \xi, H\rangle$.
\end{theo}

\pf $(a)\Rightarrow (b)$:  
Let $F_t=z-tH$ and assume that 
$\Lambda^m(P^m)=0$ for all $m\ge 1$. 
Then  by Proposition \ref{Lambda-Nil}, 
we know that $JH$ is nilpotent and,  
by Lemma \ref{H-nil-lemma}, 
$\cJ F_t(z)\equiv 1$.   

Since we have assumed that 
the  {\it Jacobian conjecture}  holds, 
the formal inverse map $G_t=z+tN_t$ of $F_t$ 
is also a polynomial map over $K[t]$. 
Then  by Proposition \ref{t-Lambda-JC}, 
we have $\Lambda^m(P^{m+1})=0$ when $m\gg 0$. 

$(b)\Rightarrow (a)$:  
Let $F$ be any polynomial map with $\cJ F\equiv 1$. 
By the well-known homogeneous reduction 
in \cite{BCW} and \cite{J}, we may assume that 
$F$ has the form $F=z-H$ such that 
$H$ is homogeneous (of degree $3$).

Under the homogeneous assumption on $H$, it is easy to 
check that the condition $\cJ F\equiv 1$ implies 
(actually is equivalent to) the statement 
that $JH$ is nilpotent. 

Then  by Proposition \ref{Lambda-Nil}, we have
$\Lambda^m(P^m)=0$ for all $m\ge 1$. 
Since we have assumed 
that the  {\it vanishing conjecture}, Conjecture \ref{VC}, 
holds for $\Lambda$ and $P(\xiz)$, we have  
$\Lambda^m(P^{m+1})=0$ when $m\gg 0$. Then  
by Proposition \ref{t-Lambda-JC}, 
we know that the formal inverse map 
$G_t(z)$ is also a polynomial map over $K[t]$. 
Hence so is $G(z)=G_{t=1}(z)$.
\epfv

\begin{rmk}
Note that  it has been shown in \cite{HNP} and \cite{GVC} that 
the  {\it Jacobian conjecture}  is actually equivalent to 
the  {\it vanishing conjecture}, Conjecture \ref{VC}, 
without any assumption on $P(\xiz)\in \cA_K[\xiz]$. 
\end{rmk}

{\small \sc Department of Mathematics, Illinois State University,
Normal, IL 61790-4520.}\, 
{\em E-mail}: wzhao@ilstu.edu.


\begin{thebibliography}{FLM2}

\bibitem[A]{Ab} S. S. Abhyankar, {\it Lectures in algebraic geometry}. Notes by Chris Christensen, Purdue University, 
1974.


\bibitem[BCW]{BCW} H. Bass, E. Connell, D. Wright, {\it The Jacobian
conjecture, reduction of degree and formal expansion of the inverse}.
Bull.  Amer. Math.  Soc.  \textbf{7}, (1982), 287--330. [MR 83k:14028]. Zbl.539.13012.


\bibitem[B1]{B} J.-E. Bj\"ork, 
{\it Rings of differential operators}. North-Holland Publishing Co., Amsterdam-New York, 1979. [MR0549189].

\bibitem[B2]{B2} 
J.-E. Bj\"ork, {\it Analytic ${\mathcal D}$-modules and applications}.  Mathematics and its Applications, 247. Kluwer Academic 
Publishers Group, Dordrecht, 1993. [MR1232191].  

\bibitem[BE]{BE} M. de Bondt and A. van den Essen, 
{\it  A reduction of the Jacobian conjecture to the Symmetric Case}.  
Proc. Amer. Math. Soc. {\bf 133} (2005), no. 8, 2201--2205. [MR2138860].

\bibitem[C]{C} S. C. Coutinho,
{\it A primer of algebraic $D$-modules}. 
London Mathematical Society Student Texts, 33. 
Cambridge University Press, Cambridge, 1995. [MR1356713].

\bibitem[D]{D} N. Budur, {\it On the $V$-filtration of $\mathcal D$-modules}. Geometric methods in algebra and number theory 59--70, Progr. Math., 235, Birkh\"auser Boston, Boston, MA, 2005. [MR2159377].

\bibitem[E]{E} A. van den Essen, {\em Polynomial automorphisms and the Jacobian conjecture}.  Progress in Mathematics, 190. Birkh\"auser Verlag, Basel, 2000. [MR1790619]. 
 
\bibitem[EWZ]{EWZ} A. van den Essen, R. Willems and W. Zhao, 
{\it Some results on the vanishing conjecture of differential operators with constant coefficients}. arXiv:0903.1478v1 [math.AC]. Under submission. 

\bibitem[EZ]{EZ} A. van den Essen and W. Zhao, 
{\it Two results on Hessian nilpotent polynomials}. J. Pure Appl. Algebra {\bf 212} (2008), no. 10, 2190--2193. [MR2418165]. See also math.arXiv:0704.1690. 

\bibitem[J]{J} A. V. Jag\v zev, {\it On a problem of O.-H. Keller.} (Russian) Sibirsk. Mat. Zh. 21 (1980), no. 5, 141--150, 191. [MR0592226].

\bibitem[K]{Ke} O. H. Keller, 
{\it Ganze gremona-transformationen}. Monats. Math. Physik {\bf 47} (1939), no.\,1, 299-306. [MR1550818].

\bibitem[M]{Me} G. Meng, {\it Legendre transform, Hessian conjecture and tree formula}.  Appl. Math. Lett. 19 (2006), no. 6, 503--510, [MR2221506]. See also math-ph/0308035.

\bibitem[S]{Sm} S. Smale, {\it Mathematical problems for the next century}. Math. Intelligencer 20, No. 2, 7-15, 1998. [MR1631413].

\bibitem[Z1]{Exp} W. Zhao, {\it Exponential formulas for the Jacobians and Jacobian matrices of analytic maps}. J. Pure Appl. Algebra {\bf 166} (2002), no. 3, 321--336. 

\bibitem[Z2]{BurgersEq} W. Zhao, {\it Inversion problem, Legendre transform and inviscid Burgers' equation},  J. Pure Appl. Algebra {\bf 199} (2005), no.\,1-3, 299--317. [MR2134306]. See also math.CV/0403020.

\bibitem[Z3]{HNP} W. Zhao, {\it Hessian nilpotent polynomials and the Jacobian conjecture}, Trans. Amer. Math. Soc. 359 (2007), no. 1, 249--274 (electronic).
[MR2247890]. See also math.CV/0409534.

\bibitem[Z4]{GVC} W. Zhao, {\it A Vanishing conjecture on differential operators with constant coefficients}, Acta Mathematica Vietnamica, vol 32 (2007), no. 3, 107--134. [MR2368014]. See also arXiv:0704.1691v2 [math.CV].

\bibitem[Z5]{IC} W. Zhao, {\it Images of commuting  differential operators of order one with constant leading coefficients}.  J. Alg. {\bf 324} (2010),  no. 2, 231--247. See also arXiv:0902.0210 [math.CV].  

\bibitem[Z6]{T-Deform} W. Zhao, {\it A deformation of commutative polynomial algebras in even numbers of variables}.  Cent. Eur. J. Math.\,{\bf 8} (2010), no.\,1, 73--97. [MR2593265]. See also arXiv:0907.3990v1 [math.AC].

\end{thebibliography}
\end{document}